\documentclass[12pt]{article}
\usepackage{amssymb}
\usepackage{latexsym}
\begin{document}
\def\Ubf#1{{\baselineskip=0pt\vtop{\hbox{$#1$}\hbox{$\sim$}}}{}}
\def\ubf#1{{\baselineskip=0pt\vtop{\hbox{$#1$}\hbox{$\scriptscriptstyle\sim$}}}{}}
\def\R{{\Bbb R}}
\def\V{{\Bbb V}}

\title{Orbit cardinals: On the {\it definable cardinalities} of quotient spaces of the 
form $X/G$, where $G$ acts on a Polish space $X$}         
\author{Greg Hjorth}        
\date{\today}          
\maketitle

{\bf $\S$0 Preface} 
This paper 
is part of a project to obtain a structure theory for the simplest and most general objects of  
mathematics. 
I wish to consider the {\it definable cardinalities} that arise from 
the continuous actions of Polish groups. The philosophy is to calculate 
cardinalities using only sets and functions that are in some sense {\it reasonably 
definable}. As with \cite{kechris3} and \cite{friedmanstanley}, the study of definable 
cardinaliites is intended to be an abstract investigation of classification 
problems, in that we may say that the classification of the equivalence relation $E$ on $X$ is 
harder than the classification of $F$ on $Y$ if the definable cardinality of 
$X/E$ exceeds that of $Y/F$.

Of course the notion of reasonably definable is vague and subject to 
personal taste and prejudice. I will choose to explicate this notion by taking perhaps the 
most generous definition in wide currency. For me, the reasonably definable sets are those 
that appear in $L({\Bbb R})$, the universe of all objects that arise from transfinite operations 
applied to ${\Bbb R}$. 

This may very well be too liberal for some, 
and an alternative approach would be to restrict ourselves to say the Borel sets, 
thereby giving us the notion of {\it Borel cardinality}; alternatively we may diet on the 
sets and functions  
arising in the $\sigma$-algebra generated by the open sets and closed under continuous images. 
For most of the problems considered below the structure one obtains for the cardinals in 
$L({\Bbb R})$ closely resembles that suggested by the Borel sets and functions. Indeed, 
under the assumption of $AD^{L({\Bbb R})}$, 
the universe of $L({\Bbb R})$ fills out the sketch outlined for us by the 
Borel sets, providing a canonical model of ZF where every set of reals has the regularity properties 
such as being Lebesgue measurable and the cardinal structure plays out the suggestions made by 
the Borel equivalence relations. 

It should be stressed that $L({\Bbb R})$ is a model of ZF, but not of choice. Thus not 
every set can be wellordered, and consequently not every cardinal corresponds to an ordinal. 
For instance, the cardinality of $2^{\aleph_0}$ is not an ordinal in $L({\Bbb R})$ -- just as there 
is no Borel wellordering of ${\Bbb R}$ in ZFC. Morever, the existence of a surjection 
$\pi:A\rightarrow B$ does not guarantee that $|A|$, the cardinality of $A$, is less than 
the cardinality of $B$, in the sense of there being an injection from $B$ to $A$. For instance, 
although there is a surjection from ${\Bbb R}$ to ${\Bbb Q}/{\Bbb R}$, there is no injection in 
$L({\Bbb R})$ from ${\Bbb Q}/{\Bbb R}$ to ${\Bbb R}$. 
To keep the distinctions in view, I will always write $|A|_{L({\Bbb R})}$ to 
indicate the cardinality of $A$ as calculated in ${L({\Bbb R})}$. 

The first result is one in a long line of generalizations of the Glimm-Effros dichotomy 
for Polish group actions.\\

0.1. Theorem($AD^{L({\Bbb R})}$). Let $G$ be a Polish group acting continuously on a Polish space 
$X$, and let $A\subset X$ be in $L({\Bbb R})$. 
Then either 

(I) $|A/G|_{L({\Bbb R})}\leq |2^{<\omega_1}|_{L({\Bbb R})}$, or 

(II) $|{\Bbb R}/{\Bbb Q}|_{L({\Bbb R})}\leq |A/G|_{L({\Bbb R})}$.\\

The proof also works in the $AD_{\Bbb R}$ context, thereby answering a question from \cite{beckerkechris}. 
By much the same argument\\

0.2. Theorem($AD^{L({\Bbb R})}$). Let $G$ be a Polish group acting continuously on a Polish space 
$X$, and let $A\subset X$ be in $L({\Bbb R})$. 
If $|A/G|_{L({\Bbb R})}\leq |
H(\kappa)^{L({\Bbb R})}|_{L({\Bbb R})}$ for some ordinal $\kappa$, then 
$|A/G|_{L({\Bbb R})}\leq |
H(\omega_1)|_{L({\Bbb R})}$.\\

So that in $L({\Bbb R})$, if there are wellorderable sets that can be assigned as complete invariants 
to the orbit equivalence relation, 
then we can in fact find elements in $HC$, the collection of all hereditarily countably sets, as 
complete invariants. One obtains this kind of classification in the Scott analysis of the isomorphism relation on 
countable structures. 

Another direction was suggested by recent work of Howard Becker's:\\

0.3. Theorem(Becker). Let $G$ be a Polish group with a left invariant complete metric 
acting continuously on a Polish space $X$. Then either 

(I) there is a Borel $\theta:X\rightarrow 2^{\omega}$ such that for all $x_1,x_2\in X$ 
\[\exists g\in G(g\cdot x_1=x_2)\Leftrightarrow \theta(x_1)=\theta(x_2)\]
or 

(II) there is a Borel $\theta: {\Bbb R}\rightarrow X$ such that for all $r_1, r_2\in {\Bbb R}$ 
\[r_1-r_2 \in {\Bbb Q} \Leftrightarrow \exists g\in G(g\cdot \theta(r_1)=\theta(r_2)).\]\\

The class of Polish groups with a left invariant complete metric includes all locally compact and 
all solvable Polish groups, but not the symmetric group of permutations on a countably infinite 
set with the topology of pointwise topology nor the automorphism group of $[0,1]$ under the 
compact-open topology. 

While Becker also established approximately this result 
for $\Ubf{\Sigma}^1_1$ sets (that is, those arising as the continuous images 
of Borel sets), he did so only under the additional assumption that {\it every real has a sharp}. 
Below we obtain just in ZFC that\\

0.4. Theorem. Let $G$ be a Polish group with a left invariant complete metric 
acting continuously on a Polish space $X$, and let $A\subset X$ be $\Ubf{\Sigma}^1_1$. Then either 

(I) there is a $\Ubf{\Delta}^1_2$ 
$\theta:A\rightarrow 2^{\omega}$ such that for all $x_1,x_2\in A$ 
\[\exists g\in G(g\cdot x_1=x_2)\Leftrightarrow \theta(x_1)=\theta(x_2)\]
or 

(II) there is a Borel $\theta: {\Bbb R}\rightarrow A$ such that for all $r_1, r_2\in {\Bbb R}$ 
\[r_1-r_2 \in {\Bbb Q} \Leftrightarrow \exists g\in G(g\cdot \theta(r_1)=\theta(r_2)).\]\\

The proof also yields under appropriate determinacy or large cardinal assumptions a generalization that 
Becker's arguments do not seem to give under any hypothesis.\\

0.5. Theorem($AD^{L({\Bbb R})}$). Let $G$ be a Polish group with a left invariant complete metric 
acting continuously on a Polish space $X$, and let $A\subset X$ be in $L({\Bbb R})$. Then either 

(I) $|A/G|_{L({\Bbb R})}\leq |2^{\omega}|_{L({\Bbb R})}$, or 

(II) $|{\Bbb R}/{\Bbb Q}|_{L({\Bbb R})}\leq |A/G|_{L({\Bbb R})}$.\\

It should be noted here that this may be viewed as a generalization of 0.3, since (I) is equivalent to 
the existence of some $\theta\in L({\Bbb R}), \theta:X\rightarrow \R$ such that for all $x_1,x_2\in A$ 
\[\exists g\in G(g\cdot x_1=x_2)\Leftrightarrow \theta(x_1)=\theta(x_2),\]
while (II) is equivalent to the existence $\theta\in L({\Bbb R})$, $\theta:\R\rightarrow X$, 
such that for all $r_1, r_2\in {\Bbb R}$ 
\[r_1-r_2 \in {\Bbb Q} \Leftrightarrow \exists g\in G(g\cdot \theta(r_1)=\theta(r_2)).\]

By the same method one obtains that there is no way the orbit structure of such a group action 
can reduce the equality relation on countable sets of reals.\\

0.6. Theorem. Let $G$ be a Polish group with a left invariant complete metric acting continuously 
on a Polish space $X$. Then there is no Borel $\theta: 2^{\Bbb N}\times\R\rightarrow X$ 
such that for all $x, y \in 2^{\Bbb N}\times\R$ 
\[\{x(n, \cdot):n\in {\Bbb N}\}=\{y(n, \cdot):n\in {\Bbb N}\}\Leftrightarrow \exists g\in G(g\cdot x=y).\]\\

In the AD$^{L(\R)}$ context this yields that 
$|{\cal P}_{\aleph_0}(\R)|_{L({\Bbb R})}\not\leq |X/G|_{L({\Bbb R})}$ -- the definable cardinality of 
the set of all countable sets of reals is not below that of the set of $G$-orbits.

Finally, since Becker's result implies Vaught's conjecture for Polish groups admitting a left invariant 
complete metric, he was led to ask whether these are the {\it only} Polish 
groups satisfying Vaught's conjecture. In answer:\\

0.7. Theorem. There is a Polish group $G$ with no compatible left invariant metric such that 
whenever it acts continuously on a Polish space $X$, either 

(I) $|X/G|\leq\aleph_0$; or 

(II) $2^{\aleph_0}\leq |X/G|$.\\

The goup arises as Aut$(M)$, for $M$ a countable model constructed by Julia Knight. 

The different sections can be read independently, with only the proofs of $\S$3 requiring a knowledge 
of determinacy. The background material is spread through $\S 1$, $\S 2$, and $\S 4$, with $\S 5$ 
requiring $\S 1$ and $\S 4$, $\S 3$ assuming $\S 2$ and $\S 1$, and $\S 5$ only $\S 1$. 

\newpage
{\bf $\S$1 On Polish groups}\\

This section collects together 
some background on Polish group actions. Further discussion, along with a few of the 
proofs and most of the references, can be 
found in \cite{beckerkechris}, \cite{kechris2}, or \cite{kechris}.\\ 

1.1. Definition. A topological group is said to be {\it Polish} if it is Polish as a topological 
space -- which is to say that it is separable and allows a complete metric. If $G$ is a Polish 
group and $X$ is a Polish space on which it acts continuously, then $X$ is said to be a {\it Polish 
$G$-space}. $E^X_G$ is the orbit equivalence on $X$, given by 
\[x_1E_G^X x_2 \Leftrightarrow \exists g\in G(g\cdot x_1 = x_2).\] 
The orbit, $G\cdot x$ of a point $x$ in $X$ is denoted by $[x]_G$. $X/G$ denotes the collection of 
orbits, $\{[x]_G:x\in X\}$.\\

1.2. Example. Let $S_{\infty}$ be the group of all permutations of the natural numbers, 
and let $2^{{\Bbb N}\times{\Bbb N}}$ be the space of all functions from 
${{\Bbb N}\times{\Bbb N}}$ to $\{0,1\}$. 
Equip $2^{{\Bbb N}\times{\Bbb N}}$ with the product topology and 
$S_{\infty}$ with the topology of pointwise convergence, under which we have that $S_{\infty}$ is 
a Polish group and $2^{{\Bbb N}\times{\Bbb N}}$ is a Polish $S_{\infty}$-space under the 
action defined by 
\[(g\cdot x)(n,m)=x(g^{-1}(n),g^{-1}(m))\]
for any $x\in 2^{{\Bbb N}\times{\Bbb N}}$ and $g\in S_{\infty}$. 

There is a natural sense in which we may view elements of $2^{{\Bbb N}\times{\Bbb N}}$ as coding 
countable structures whose underlying set is ${\Bbb N}$ and whose only relation is a single 
binary relation, the extension of which equals $\{(m,n):x(m,n)=1\}$. If for $x\in 2^{{\Bbb N}\times{\Bbb N}}$ 
we let ${\cal M}_x$ be the corresponding model then we obtain that for all $x_1, x_2 $ in the space 
\[{\cal M}_{x_1}\cong {\cal M}_{x_2}\Leftrightarrow \exists g\in S_{\infty}(g\cdot x_1=x_2).\]
This can be extended in a simple minded fashion to allow elements of $2^{({\Bbb N}^{<{\Bbb N}})}$ 
to code models of an arbitrary countable language, and to let $S_{\infty}$ act so that it again 
induces isomorphism as its orbit equivalence relation.\\

In analyzing $S_{\infty}$ it is often possible to use model theoretic ideas, such as 
types; in the context of arbitrary Polish group actions we can hope instead to use 
the notion of {\it Vaught transforms}.\\

1.3. Definition. Let $G$ be a Polish group and $X$ a Polish $G$-space. Then for $B\subset X$, 
$U\subset G$ open, $B^{\Delta U}$ is the set of $x\in X$ such that for a non-meager set of 
$g\in U$, $g\cdot x \in B$; $B^{* U}$ is the set of $x\in X$ such that for a comeager set of 
$g\in U$, $g\cdot x \in B$.
For $x\in X$, $B\subset X$, $U\subset G$ open, $\forall^*g\in U(g\cdot x\in B)$ indicates that 
for a relatively comeager set of $g\in U$, $g\cdot x\in B$.  $\exists^*g\in U(g\cdot x\in B)$ is used 
to indicate that for a non-meager set of $g\in U$, $g\cdot x\in B$.\\

It is generally only advisable to consider the Vaught transform of $B$ when it is sufficiently 
well behaved to guarantee that the transforms have the Baire property -- for instance, if $B$ is 
Borel, or in $L({\Bbb R})$ under suitable hypotheses. 
In the case that, say, $(B_i)$ is a sequence of Borel sets
\[(\bigcup B_i)^{\Delta U}=\bigcup\{(B_i)^{*V}:V\subset U, V\neq \emptyset, i\in{\Bbb N}\},\]
and thus we obtain that the Vaught transform of a Borel set is again Borel. 

For general equivalence relations, induced by a group actions, or arising in some other manner, 
there is a spectrum of ways in which they may be compared, of which I mention those that 
will be most important in the remainder of the paper.\\

1.4. Definition. For $E$ and $F$ equivalence relations on Polish spaces $X$and $Y$, 
$E\leq_B F$, $E$ {\it Borel reduces} $F$, indicates that 
there is a Borel function $\theta:X\rightarrow Y$ such that 
for all $x_1, x_2\in X$
\[x_1E x_2\Leftrightarrow \theta(x_1)F\theta(x_2);\]
we write $E\leq F_c$, $E\leq_{\Ubf{\Delta}^1_2} F$, $E\leq_{L({\Bbb R})} F$ to 
indicate that there is, respectively, a continuous, $\Ubf{\Delta}^1_2$, or 
$L({\Bbb R})$, $\theta:X\rightarrow Y$ such that 
for all $x_1, x_2\in X$
\[x_1E x_2\Leftrightarrow \theta(x_1)F\theta(x_2);\]
here that we may assume without loss of generality that $X$ and $Y$ are in $L({\Bbb R})$, 
the smallest class inner model of ZF containing the reals. 
One writes $E\sqsubseteq_B F$, $E\sqsubseteq_c F$, and $E\sqsubseteq_{L({\Bbb R})} F$  
if there is a one-to-one $\theta$ that 
performs the above described reduction, and is Borel, continuous, 
or in $L({\Bbb R})$ respectively. These notions are graded, since all continuous functions 
are Borel, all Borel are ${\Ubf{\Delta}^1_2}$, which in turn lie inside $L({\Bbb R})$.\\

In this paper I will only be interested in the reductions above. These suggest a 
notion of equivalence among equivalence relations, defined to hold 
when we have bi-reducibility, in the sense that $E$ and $F$ are 
Borel equivalence and $E\leq_B F\leq_B E$. We might also define a rival notion of equivalence 
to hold when there is a Borel bijection $\theta$ between the underlying Borel 
spaces $X$ and $Y$ with 
\[\forall x_1,x_2\in X(x_1E x_2\Leftrightarrow \theta(x_1)F\theta(x_2)).\]
It turns out that the definition at 1.4 better reflects the idea of {\it definable cardinality}.\\

1.5. Definition. $E_0$ is the equivalence relation of eventual agreement on infinite 
sequences of $0$'s and $1$'s, so that for $x,y\in 2^{\Bbb N}$ 
\[xE_0 y \Leftrightarrow \exists N\forall n>N(x(n)=y(n)).\]
It is known that under the ordering of Borel reducibility, $E_0$ is equivalent to the 
more familiar Vitali equivalence relation given by 
\[x E_v y\Leftrightarrow |x-y|\in {\Bbb Q},\]
in as much as $E_0\leq_BE_v\leq_B E_0$. 

For the first inequality it suffices to consider $(0,1)/{\Bbb Q}$. Then let $p_n$ denote the 
$n^{th}$ prime, and for $x\in (0,1)$ let $\theta(x)=(r_n(x))_{n\in{\Bbb N}}$ denote the decimal 
expansion of $x$ with respect to the varying basis $(p_n)_{n\in{\Bbb N}}$ -- so that 
$\Sigma\{r_n(x)/p_n\cdot p_{n-1}\cdot...p_0:n\in{\Bbb N}\}=x$ and each $r_n(x)\in\{0,1,2,...,p_n-1\}$; 
in the case of there being more than one such expansion -- which corresponds to a recurring $9$ in an 
infinte decimal expansion - we can convene to choose the expansion that 
terminates with $r_n(x)=0$ for all sufficiently large $n$. 
By the uniqueness of the $(p_n(x))$ with the above properties, 
we obtain that $x_1-x_2\in{\Bbb Q}$ if and 
only if $\theta(x_1)$ and $\theta(x_2)$ eventually agree. From here we can organize a coding by 
elements in $2^{\omega}$, with similar properties and hence a reduction to $E_0$. 
(I am very grateful to Itay Neeman for pointing out this short proof.) 

For the second inequality, let $(q_n)_{n\in\omega}$ list the rationals, and choose a family 
$(V_s)_{s\in 2^{<\omega}}$ of open sets such that $s\subset t\Rightarrow V_s \supset \overline{V}_t$, 
$s(n)\neq t(n) \Rightarrow q_n\cdot V_s\cap V_t=0$, and 
for $lh(s)=lh(t)=n$ $w\in 2^{<\omega}$ $q_n\cdot V_{sw}=V_{tw}$, where $sw$ refers to the concatenation 
of $s$ followed by $w$. Then $\theta$ with by $\{\theta(x)\}=_{df}\bigcap V_{x|n}$ for $x\in 2^{\omega}$ 
defines the reduction.

While from the point of ZFC cardinals, $2^{\Bbb N}/E_0$ (or ${\Bbb R}/E_v$) both 
have cardinality $2^{\aleph_0}$, and hence the same size as 
$2^{\Bbb N}$, or ${\Bbb R}$, from the point of view of {\it definable cardinals}, 
these sets are very different. For instance, in 
$L({\Bbb R})$ there is no injection from $2^{\Bbb N}/E_0$ to $2^{\Bbb N}$. 
Similarly from the context of Borel structure, there is no Borel $\theta:2^{\Bbb N}\rightarrow 
2^{\Bbb N}$ such that for all $x_1,x_2$ 
\[x_1 E_0 x_2 \Leftrightarrow \theta(x_1)=\theta(x_2).\]
Here id$(2^{\omega})$ is the equality relation on $2^{\omega}$, which, as the 
collection of sequences from $\{0,1\}$, may be identified with 
$2^{\Bbb N}$. 
id$(2^{<\omega_1})$ is the equality relation on countable transfinite sequences 
of $0$'s and $1$'s. 
Again, while $2^{\omega}$ and $2^{<\omega_1}$ have the same cardinality in ZFC, there is 
no {\it reasonably definable} injection from $2^{\omega}$ to $2^{<\omega_1}$; so, under 
suitable large cardinal assumptions, we have, for instance, no such injection in 
$L({\Bbb R})$.\\

1.6. Definition. HC denotes the collection of all sets whose transitive cardinality is 
countable -- that is, every $x_0\in x$ is countable, every $x_1\in x_0\in x$ is countable, 
and so on.\\

It is known from the Scott analysis of \cite{keisler} and the more recent results 
of \cite{beckerkechris} that if $X$ is a Polish $S_{\infty}$-space then $E_{S_{\infty}}^X\leq
_{\Ubf{\Delta}^1_2}$id(HC), in the sense of there being an ${\Ubf{\Delta}^1_2}$ {\it in the 
codes} function; that is to say, there is a ${\Ubf{\Delta}^1_2}$ function from $X$ to 
elements of $2^{{\Bbb N}\times{\Bbb N}}$, such that any two $x_1, x_2\in X$ are orbit 
equivalent if and only $\theta(x_1)$ and $\theta(x_2)$ code the same element in HC. 

Another theorem along these lines is due to Hjorth-Kechris and, independently, Becker. This result 
concerns the kind of classification one finds with the Ulm invariants for 
countable abelain $p$-groups, and states 
that if $G$ is any Polish group and $X$ is a Polish $G$-space, then either $E_0\sqsubseteq_c E^X_G$ 
or $E^X_G\leq _{\Ubf{\Delta}^1_2}$id$(2^{<\omega_1})$; a proof can be found in \cite{hjorthkechris}.\\

1.7. Definition. A Polish group $G$ is said to be a {\it cli group} if it has a compatible 
left invariant complete metric -- that is to say there is a compatible complete metric $d$ such 
that for all $g, h_1, h_2\in G$ 
\[d(h_1,h_2)=d(gh_1, gh_2).\]\\

It is known that all Polish groups have a compatible left invariant metric, but not all have 
a {\it complete} left invariant metric. For instance neither $S_{\infty}$ nor the homeomorphism 
group of the unit interval are cli groups. On the other hand, all abelian and locally compact groups 
are cli groups. A group has a left invariant complete matric if and only if it has a right invariant 
complete metric, since we can pass from one to the other by setting $d^*(g,h)=d(g^{-1},h^{-1})$. 

For left invariant metrics the notion convergence is topological: $(g_i)_{i\in{\Bbb N}}$ will 
Cauchy if and only if for each open neighbourhood $U$ of the identity there is some $N$ 
such that for all $n, m\geq N$ $g_n^{-1}g_m\in U$. Thus in particular, if one left invariant metric is 
complete they all are. 

The following important fact appears in \cite{sami}:\\

1.8. Lemma. Let $G$ be a Polish group and $X$ a Polish $G$-space. Then for $x\in X$ we let 
 the {\it stabilizer of $x$}, $\{g\in G:g\cdot x=x\}$ be denoted by 
$G_x$. Then $[x]_G$ is uniformly Borel in 
$x$ and any real coding $G_x$.\\

Finally:\\

1.9. Theorem(Effros). $X$ and $G$ as in 1.8, $x\in X$. $[x]_G\in \Ubf{\Pi}^0_2$ if and only if 
the map $G\rightarrow [x]_G$, $g\mapsto g\cdot x$ is open.

\newpage
{\bf $\S$2 In $L({\Bbb R})$}\\

We will need much the same technology as employed in \cite{hjorth}, 
but working with arbitrary Polish spaces. Here I will assume that 
given a point $x$ in 
some Polish space $X$ the reader is willing to allow that we 
can make sense of constructing from $x$ and forming 
the smallest inner class model, $L[x]$, containing $x$. Strictly speaking we need instead 
to fix a real $z$ coding a presentation of $X$, and speak of constructing from the 
pair $(z, y(x))$, where $y(x)$ is an element of $2^{\omega}$ that codes $x$ relative to 
the presentation given to us by $z$. Instead of being precise and strict, I will 
be more informal and treat the elements of any Polish space in exactly the same fashion 
as the {\it recursive Polish spaces}, such as ${\Bbb R}$ and 
$2^{\omega}$ -- keeping in the background that this is not quite accurate but more concise 
and easily rectifiable. Alternatively, the reader may interpret the results below as holding 
only for the recursive Polish spaces,  in the sense of \cite{moschovakis}, 
but allowing the usual relativizations to a parameter. 
Finally, we may use the fact that all uncountable Polish spaces are Borel isomorphic to code any 
Polish space by elements of $2^{\Bbb N}$, again with the coding taking place relative to some 
parameter $z$.

However the reader chooses to explicate the notion of constructing from a point in a Polish space or 
using a parameter to code such a space, the notation $x\in X\cap M$, for $M$ an inner 
model, means that $x$ is a point coded 
by a real in $M$, and that the parameter used to code $X$ exists in $M$. We can think of $U$ as being 
an open set coded in $M$ if there are sequences $(q_i)_{i\in{\Bbb N}}$ of rationals 
$(x_i)_{i\in {\Bbb N}}$ of points in $X$, both in $M$, such that $U$ equals the set of elements $x\in X$ for  
which there is some $i$ with $x$ within distance $q_i$ of $x_i$.

Up to isomorphism, all Polish spaces exist in $L(\R)$. Thus it will be convenient to have a standing 
assumption that all our Polish space are in fact an elements of this 
inner model; this assumption can always be made without loss of generality. 

The theory of $L(\R)$ will be developed under the determinacy assumption AD$^{L(\R)}$ which states 
that every subset of $\omega^{\omega}$ in $L(\R)$ is determined -- 
one of players has a winning strategy in 
the infinite game where I and II alternate in playing integers, and the victor is decided on the 
basis of whether the resulting element of $\omega^{\omega}$ is in the specified subset. 
While ZFC alone is too weak to decide most of the natural questions regarding $L(\R)$, 
the assumption of AD$^{L(\R)}$ provides a canonical theory for this inner model. 
There is a widespread acceptance of this assumption in the study of $L(\R)$ among 
set theorists -- partly because it leads to a theory for the sets of reals in $L(\R)$ 
which continues pattern we find for the Borel sets under ZFC, and partly because 
AD$^{L(\R)}$ was shown in \cite{woodin} to follow from 
large cardinal assumptions, such as the existence of a 
supercompact.\\

2.1. Definition. If $X$ is a Polish space, 
$A\subset X$ is said to be $\infty$-Borel if there is an ordinal $\alpha$, a set $S\subset \alpha$, 
and a formula $\varphi$ such that $A$ equals
\[\{x\in X: L_{\alpha}[x, S]\models \varphi(x, S)\}.\]\\

2.2. Theorem(Woodin). Assume $AD^{L({\Bbb R})}$ and let $A\subset X$ be OD$^{L({\Bbb R})}_{x}$ for 
some real $x\in {\Bbb R}$. Then there is an $\infty$-Borel code for $A$ in HOD$^{L
({\Bbb R})}_{x}$. (See \cite{schimmerling}.)\\

2.3. Definition. Let $C$ be a transitive set in $M$, a class inner model. 
Then HOD$^{M}_C$ is the smallest 
inner class model containing $C$ and closed under ordinal definability (as calculated from the 
point of view of $M$). OD$^{M}_C$ denotes all sets that are definable over $M$ in the usual Levy hierarchy 
from an ordinal, $C$, and finitely many elements of $C$. 
For $X$ a Polish space whose presentation exists in 
$M$, ${\Bbb B}(C,X,M)$ denotes $\{A\in ({\cal P}(X))^M: A\in{\rm OD}^{M}_C\}$, and can be 
viewed as a Boolean algebra in the natural sense. 
(Note that $({\cal P}(X))^M$ refers not to the true power set of $X$, but the power set of 
$X\cap M$ inside $M$.) For $G\subset {\Bbb B}(C,X,M)$ a sufficiently generic filter, 
we may define a point $x(G)\in X$ by the requirement that for all open $U\subset X$ whose 
code exists in $M$ we have that 
\[x\in U\Leftrightarrow U\in G.\]\\

The statement of the next theorem is slightly more general than is usual;  
the proof however follows exactly as does the usual proof, given in \cite{hjorth}. 
The one variation is that here our inner model HOD$^M_C$ need not satisfy choice.\\

2.4. Theorem(Vopenka). Fix $M$, $C$, and $X$ as above, and assume that $C$ includes a code for 
$X$. Then there exists ${\Bbb B}$ in 
HOD$_C^M$, $i:{\Bbb B}\cong{\Bbb B}(C, X, M)$ in $M$, such that: 

(i) for all $x\in X\cap M$, $G(x)=_{df}\{i^{-1}(A): x\in A, A\in {\rm OD}^M_C\}$ is 
HOD$_C^M$-generic for ${\Bbb B}$;

(ii) there is a HOD$_C^M$-generic for ${\Bbb B}$ below every non-zero element in $M$; 

(iii) if $H\subset {\Bbb B}$ is HOD$_C^M$-generic, if we let $G=i[H]$ then 
$x(G)\in$ HOD$^M_C[H]$, and for all ordinals $\alpha$, $\vec c \in C$, and formulas 
$\varphi$, 
\[L_{\alpha}(C, x(G))\models \varphi(\vec c, x(G))\Leftrightarrow 
\{x\in X\cap M: \varphi(\vec c, x)\in G\};\]

(iv) $i\in$OD$^M_C$.\\

Note that ${\Bbb B}$ from this theorem will have size at most ${\cal P}(X)^M$ -- the 
set of all subsets of $X$ in $M$ -- and so will be have size at most $(2^{2^{\aleph_0}})^M$. 
The following important result may be found in \cite{martinsteel}:\\

2.5. Theorem(Martin-Moschovakis-Steel). Assume $AD^{L({\Bbb R})}$. Then in 
$L({\Bbb R})$, Scale($\Ubf{\Sigma}^2_1)$. (See \cite{martinsteel}.)\\

It follows from entirely general facts that every non-empty 
$\Sigma^2_1$ collection of 
sets of reals in $L({\Bbb R})$ has a $\Ubf{\Sigma}^2_1$ member, and thus 
there will be a member of this collection which is the projection of a {\it tree}, in the sense of 
\cite{moschovakis}.\\

2.6. Definition. For $A\in L({\Bbb R})$, $|A|_{L({\Bbb R})}$ denotes the 
cardinality of $A$ as calculated in $L({\Bbb R})$. So $|A|_{L({\Bbb R})}
\leq|B|_{L({\Bbb R})}$ if there is an injection from $A$ to $B$;  
by Schroeder-Bernstein, they have the same cardinality  only if there is a bijection 
between them. For $\kappa$ an ordinal, 
H$(\kappa)$ denotes the collection of sets whose transitive closure has 
size less than $\kappa$. Thus the class of all wellorderable sets in $L({\Bbb R})$ is the union 
$\bigcup_{\kappa\in Ord}($H$(\kappa)$)$^{L({\Bbb R})}$. HC equals H$(\omega_1)$.\\

Here it is worth collecting together some facts about $L({\Bbb R})$-cardinals under 
the assumption of AD$^{L({\Bbb R})}$. Note that $H(\omega_1)=H(\omega_1)^{L(\R)}.$\\

2.7. Theorem(folklore). Assume $AD^{L({\Bbb R})}$. Then 

(i)  $|{\Bbb R}|_{L({\Bbb R})}\leq
|2^{\omega}|_{L({\Bbb R})}\leq |{\Bbb R}|_{L({\Bbb R})}$;

(ii) $|{\Bbb R}|_{L({\Bbb R})}\not\leq |\omega_1|_{L({\Bbb R})}\not\leq |{\Bbb R}|_{L({\Bbb R})}$; 

(iii) $|{\Bbb R}|_{L({\Bbb R})}<|{\Bbb R}/{\Bbb Q}|_{L({\Bbb R})}$; 

(iv) $|{\Bbb R}/{\Bbb Q}|_{L({\Bbb R})}\not\leq |(2^{\alpha})^{L(\R)}|_{L(\R)}$ for any ordinal 
$\alpha$; 

(v) $|{\Bbb R}/{\Bbb Q}|_{L({\Bbb R})}\leq |2^{\omega}/E_0|_{L(\R)}\leq|{\Bbb R}/{\Bbb Q}|_{L({\Bbb R})};$ 

(vi) $|\omega_1|_{L(\R)}<|2^{<\omega_1}|_{L(\R)}<|H(\omega_1)|_{L(\R)}$; 

(vii) $|2^{\omega}|_{L(\R)}<|2^{<\omega_1}|_{L(\R)}$.

Proofs. (i): This is clear even without any sort of determinacy assumptions, since there are 
Borel injections both ways. 

(ii): As can be found in \cite{jech}, there is a 
countably complete ultrafilter on $\omega_1$ under AD, so there can 
be no $\omega_1$ sequence of reals in $L(\R)$. 

(iii): $|{\Bbb R}|_{L({\Bbb R})}\leq|{\Bbb R}/{\Bbb Q}|_{L({\Bbb R})}$ since we can find a map 
from $\R$ to $\R$ such that any two distinct reals have images that are mutually generic over 
$L_{\omega_1^{ck}}$. To see the failure of reducibility in the other direction note that by 
the Lebesgue density theorem any ${\Bbb Q}$-invariant Lebesgue measurable function from 
$\R$ must be constant almost everywhere; since all functions are Lebesgue measurable in 
$L(\R)$, this suffices. 

(iv): Let $\theta:\R\rightarrow 2^{\alpha}$ be ${\Bbb Q}$-invariant and in $L(\R)$. 
Then, as in the proof of (iii), for each $\beta$ less than $\alpha$, the set 
$\{x\in \R:\theta(x)=1\}$ is either null or co-null. By Fubini's theorem in 
$L(\R)$ and all sets Lebesgue measurable, wellordered intersections of co-null sets are 
co-null, and so $\theta$ must be constant almost everywhere. 

(v): This follows as in the remarks after 1.5, since we have $E_v\leq_B E_0\leq_B E_v$. 

%

(vi): Any countable ordinal $\alpha$ can be coded be a function in $2^{<\omega_1}$ that has domain 
$\alpha$ and takes constant value $1$. The other inequality follows since 
$2^{<\omega_1}\subset H(\omega_1)$.

(vii): The non-reduction follows since there is no $\omega_1$ sequence of reals.$\Box$\\

A number of results similar to (iii), (vi), and (vii) are presented in \cite{friedmanstanley}.\\

2.8. Theorem(Woodin). Assume $AD^{L({\Bbb R})}$. If $E\in L({\Bbb R})$ on ${\Bbb R}$, then 
exactly one of the following hold: 

(I) id$(2^{\omega})\sqsubseteq_c E$; or,

(II) for some ordinal $\kappa$, $E\leq_{L(\R)}$id$(\kappa)$.\\

2.9 Corollary to the proof. Assume $AD^{L({\Bbb R})}$. Then for any set $A$, exactly one of the following holds: 

(I) $|\R|_{L(\R)}\leq |A|_{L(\R)}$; or, 

(II) for some ordinal $\kappa$, $|A|_{L(\R)}\leq |\kappa|_{L(\R)}.$\\

2.10. Theorem(Hjorth). Assume $AD^{L({\Bbb R})}$. If $E\in L({\Bbb R})$ on ${\Bbb R}$, then 
exactly one of the following hold: 

(I) $E_0
\sqsubseteq_c E$; or,

(II) for some ordinal $\kappa$, $E\leq_{L(\R)}$id$(2^{\kappa})$. (See \cite{hjorth}.)\\

2.11. Corollary to the proof. Assume $AD^{L({\Bbb R})}$. Then for any set $A$, exactly one of the following holds: 

(I) $|\R/ E_0|_{L(\R)}\leq |A|_{L(\R)}$; or, 

(II) for some ordinal $\kappa$, $|A|_{L(\R)}\leq |2^{\kappa}|_{L(\R)}.$\\

Theorem 2.10 follows by arguments similar to those used in proving 2.8. 
It is unkown whether there is 
analogue of these result for the cardinality of 
HC in $L(\R)$, but it is known that any such result would 
need to be considerably more complex.\\

2.12. Lemma. Assume $AD^{L({\Bbb R})}$. Let $E$ and $F$ be Borel -- or even $\Ubf{\Delta}^2_1$, 
or even projective --  
equivalence relations on Polish spaces $X$ and $Y$. 
Then $E\leq _{L(\R)} F$ if and only if $|X/E|_{L(\R)}\leq |Y/F|_{L(\R)}$. 

Proof. The {\it only if} direction is immediate,  so suppose that $|X/E|_{L(\R)}\leq |Y/F|_{L(\R)}$. 
Then we can find a set $R\subset X\times Y$ in $L(\R)$ so that: 

(i) $\forall (x_1, y_1), (x_2, y_2)\in R$, $x_1E x_2$ if and only if $y_1 Fy_2$; 

(ii) $\forall x\in X\exists y\in Y((x,y)\in R)$. 

Thus by 2.5 we can find such a set $R$ with $R\in \Ubf{\Sigma}^2_1$, and then a tree 
$T$ on some ordinal with $p[T]=R$. Then by the absoluteness of illfoundedness for 
trees, we can find in each model $L[T,x]$ some $y\in Y$ with 
\[(x,y)\in p[T].\]
Note here that we can define from $x$ and $T$ a wellorder of $L[T,x]$. 
Thus we may define $\theta: X\rightarrow Y$ by letting 
$\theta(x)$ be the 
first $y$ above in the canonical wellorder of $L[T,x]$.$\Box$\\

Consequently it is natural to use the ordering $\leq_{L(\R)}$ -- and by analogy 
$\leq_B$ -- in comparing Borel equivalence relations, since this is the notion of 
comparison that corresponds to definable cardinality. 

On the other hand:\\

2.13. Lemma(folklore). Let $x\in L(\R)$ be a non-empty set. Then there is a 
$\pi \in L(\R)$ and ordinal $\alpha$ such that 

(i) $\pi:\R\times \alpha \rightarrow X$ is onto; and thus 

(ii) there is a sequence $(E_{\beta})_{\beta\in\alpha}$ of equivalence relations in $L(\R)$ 
and $A\subset \{([x]_{E_{\beta}},\beta):\beta<\alpha,x\in \R\}$, and 
a bijection $\sigma:A\rightarrow X$, $A, \sigma\in L(\R)$.\\

And therefore the study of cardinalities in 
$L(\R)$ is largely the study of definable equivalence relations 
and their corresponding quotient spaces. 

The following result, stated in a rather narrow form, places the results from 
section 3 in context.\\

2.14. Theorem(Becker-Kechris). Assume AD$^{L(\R)}$. Let $G$ be a Polish group acting 
continuously on a separable metric space $X$. Then $E^A_G\sqsubseteq E^X_G$ for some 
Polish $G$-space $X$.(See \cite{beckerkechris}.)\\

\newpage

{\bf $\S$3. Generalized Ulm-type dichotomies}\\

The next two theorems are stated under entirely abstract hypotheses, assuming ZF, DC -- the 
axiom of dependent choice -- and some manner of exotic regularity property for the relevant 
sets of reals. Of course, the main interest is in the consequences for $L(\R)$, and the 
precise statements below are of technical interest.\\

3.1. Theorem. Assume ZF, DC, all sets of reals are $\infty$-Borel, and that there is 
no $\omega_1$ sequence of reals. Let $G$ be a Polish group, $X$ be a Polish $G$-space, 
and $A\subset X$ $G$-invariant. 
Then either: 

(I) $E^X_G|_A\leq$id$(2^{<\omega_1})$; or, 

(II) $E_0\sqsubseteq_c E^X_G|A$. 

Remark: Here the unadorned $\leq$ means that there just outright exists a reduction $\theta$, with 
no special definability assumption. In the context of ZF+$\neg$AC this notion has content. 

Proof.  
Let $\varphi$, $\alpha$, and $S\subset\alpha$ witness the definition of 
$\infty$-Borel, so that $A$ is equal to the set of $x\in X$ such that 
\[L_{\alpha}[S,x].\] 
Without loss of generality $S$ codes $X$, $G$, and the action, in some appropriate sense.
Since there is no $\omega_1$-sequence of reals, $\omega_1$ is strongly inacessible in $L[S,x]$ 
for any $x\in X$. Thus in particular almost every $g\in G$ is generic over $L[S,x]$ for 
the forcing notion that uses the non-empty basic open sets of $G$ ordered under inclusion as a 
forcing notion. This notion is equivalent to Cohen forcing and homogenous; thus 
as in the standard development of forcing, presented by \cite{jech}, the corresponding HOD 
of the generic extension is decided in the ground model, and 
\[ \forall ^*g\in G\forall ^*h\in G({\rm HOD}_S^{L[S, g\cdot x]}={\rm HOD}_S^{L[S, h\cdot x]}).\]
Let $M_S^x$ denote this common model, so that $\forall ^*g\in G(M^x_S={\rm HOD}_S^{L[S, g\cdot x]})$. 
Note then that $M^x_S$ depends only on $[x]_G$. 
Let $\gamma(x)$ be such that for a comeager set of $g\in G$, 
$(\beth_{\omega+\omega})^{L[S,g\cdot x]}=\gamma(x)$. Using that $M^x_S$ has a uniformly 
$\Sigma_2(S)$ wellorder, 
we can find $\theta_x\in 2^{<\omega_1}$ coding the model ${M_S^x}$ up to $\gamma(x)$ -- so that for  $(\varphi_n)_{n\in\omega}$ some reasonable enumeration of 
the formulas of set theory, $\theta_x=\{(n,\vec \alpha):L[S, g\cdot x]\models \varphi_n(S,\vec \alpha), 
\vec \alpha<\gamma(x)\}$ for a comeager set of $g\in G$. 
Note that $x\mapsto \theta_x$ is $G$-invariant. 

For $x\in A$ let ${\Bbb B}_S^x \in M^x_S$ be as 
indicated by 2.4, so that ${\Bbb B}_S^x\cong {\Bbb B}(S, X, L[S, g\cdot x])$ for a 
comeager set of $g\in G$ and is the first such algebra in the canonical wellorder of 
$M^x_S$. Then let $B_x$ be the set of $y\in X$ such that for some $g\in G$, $p\in {\Bbb B}_S^x $ 
with $g\cdot y$ $M^x_S$-generic over ${\Bbb B}_S^x $ below $p$ with 
$p$ forcing that $x(\dot{G})\in A$, 
in the sense that  
\[p\Vdash L_{\alpha}[S,x(\dot{G})]\models \varphi(S,x(\dot{G})),\]
where $\dot{G}$ refers to the generic object, and as at 2.3 $x(\dot{G})$ refers to 
the generic real. 
Let $\overline{\theta}_x=(\theta_x,p_x)$ where $p_x$ is maximal such that 
$p_x\Vdash L_{\alpha}[S,x(\dot{G})]\models \varphi(S,x(\dot{G}));$ 
$p_x$ exists by completeness of the Boolean algebra over $M^x_S$. 
$B_x$ is uniformly $\Delta^1_1(\overline(\theta)_x)$ -- in that it is 
uniformly $\Delta^1_1(w)$ for any real $w$ that codes $\overline{\theta}_x$. 
Note also that $x\mapsto B_x$ is $G$-invariant and each $B_x\subset A$. 

Now if for any $x$ we have $E_0\sqsubseteq_c E_G^X|_{B_x}$ then certainly $E_0\sqsubseteq_c A$, 
and the proof is finished. 

So suppose otherwise. 

Note that for any $y\in B_x$, $\forall ^*g\in G (\omega_1^{L[S, g\cdot y]}<\gamma(x)<\omega_1)$, and 
thus the stabilizer of $g\cdot y$, $G_{g \cdot y}=\{h\in G: hg\cdot y=g\cdot y\}$ is uniformly 
$\Delta^1_1(\overline{\theta}_x,g\cdot y)$ since it is a $\Pi^1_1(g\cdot y)$ 
singleton. Now note that $G_y=g^{-1}G_{g\cdot y} g$ will 
be a $\Sigma^1_1(\overline{\theta_x}, g\cdot y)$ singleton 
for any such $g$, and hence it will be a $\Sigma^1_1(\overline{\theta_x})$ singleton, 
and hence $\Delta^1_1(\overline{\theta_x})$. 

Thus, by 1.8, we have that $E_G^X|_{B_x}$ is uniformly Borel in any code for 
$\overline{\theta}_x$. Since $E_0\not\sqsubseteq_c E_G^X|_{B_x}$, by 
\cite{hakelo} there is a $\Delta^1_1(w)$ seperating family for any $w$ coding 
$\overline{\theta}_x$ -- in the sense that there is $(W_n)_{n\in\omega}$ a family of $G$-invariant 
$\Delta^1_1(w)$ sets that for any $y_0,y_1\in B_x$
\[y_0 E^X_G y_1\Leftrightarrow \forall n(y_0\in W_n\Leftrightarrow y_1\in W_n).\]

Thus as in \cite{hjorthkechris}, for any $y\in B_x$ we can let $\theta_0(y)$ be the set 
of terms $\tau\in M_x^S$ in ${\Bbb P}=$ Coll$(\omega,\overline{\theta}_x)$ and $q\in {\Bbb P}$ 
such that $q$ forces that $\tau$ is a $G$-invariant Borel set containing $y$. 
It is routine to use $\overline{\theta}_x$ to encode $\theta_0(y)$ as a bounded subset of 
$\omega_1$. 

Thus we may at last find $G$-invariant $\theta_1:A\rightarrow 2^{<\omega_1}$ such 
that for any $x_1, x_2$, 
\[\theta_1(x_1)=\theta_1(x_2)\Leftrightarrow
\overline{\theta}_{x_1}=\overline{\theta}{x_2}\wedge\theta_0(x_1)=
\theta_0(x_2).\]
Thus by the properties established along route, 
\[\theta_1(x_1)=\theta_1(x_2)\Leftrightarrow x_1E^X_G x_2\]
for all $x_1, x_2\in A$.$\Box$\\

3.2. Corollary. Assume AD$^{L(\R)}$. Let $G$ be a Polish group and $X$ a Polish 
$G$-space, and let $A\subset X$ be in $L(\R)$. Then either:

(I) $E^X_G|_A\leq_{L(\R)}$id$(2^{<\omega_1})$; or 

(II) $E_0\sqsubseteq_cE^X_G|_A$. 

Proof. By 2.2 we have all sets of reals $\infty$-Borel in $L(\R)$. 
As at 2.7, the determinacy assumption implies 
there is no $\omega_1$-sequence of reals, 
and we have the assumptions of 3.4.$\Box$.\\

In unpublished work Woodin has previously shown AD$_{\R}$ implies 
all sets are the projection of some tree on some 
ordinal, and hence $\infty$-Borel; thus we obtain a proof that under ZF+DC+AD$_{\R}$ either 
$E^X_G|_A\leq$id$(2^{<\omega_1})$, or $E\sqsubseteq_cE^X_G|_A$, thereby answering  
question 8.1.2 from \cite{beckerkechris}. Alternatively:\\

3.3. Theorem. Assume ZF+DC+"all sets of reals have Baire property." Let $G$ be a Polish group and let
$X$ be a Polish 
$G$-space. Let $A\subset X$ be $G$-invariant with $A=p[T]$ for $T$ a tree on some ordinal $\kappa$. 
Then either: 

(I) $E^X_G|_A\leq$id$(2^{<\omega_1})$; or, 

(II) $E_0\sqsubseteq_c E^X_G|_A$.

Proof. Let $z$ be a parameter coding the group and the action. 
Note that for $x\in A$ there will be a leftmost branch witnessing this, 
in the sense that there will $f$ with $(x,f)\in[T]$ and such that for all other 
$f_0\neq f$, either $(x,f_0)$ is not in $[T]$, or there is some $n\in\omega$ 
with $f|_n=f_0|_n$ but $f(n)<f_0(n)$. 

For $x\in A$, let $\sigma(x)$ be the 
set of nodes $s\in T$ such that $\exists ^*g\in G$(the leftmost branch witnessing $g\cdot x\in p[T]$ has 
$f\supseteq s$). By all sets having Baire property, this equals the collection of $s\in T$ such that 
for some $U\subset G$ open and non-empty, $\forall ^*g\in U$(the leftmost branch witnessing 
$g\cdot x\in p[T]$ has 
$f\supseteq s$). $x\mapsto\sigma(x)$ is $G$-invariant by the nature of the definition and the fact that 
the notions of categoricity in $G$ are invariant under translation by any $g\in G$. 
$\sigma(x)$ is countable for any $x$ by the chain condition on the ideal of meager sets in $G$. 

Thus $\sigma(x)$ may be coded by a countable tree $T_x\cong\sigma(x)$ on some ordinal $\alpha_x<\omega_1$ -- 
again with $x\mapsto T_x$ $G$-invariant. Note now that by all sets with Baire property 
$\forall^*g\in G$, the leftmost branch in 
$T$ witnessing $g\cdot x \in A$ must be in $\sigma(x)$, and thus 
$g\cdot x\in p[T_x]$. 

Now choose $\beta_x$ to be least so that 
$\forall^*g \in G((|T_x|^+)^{L[g\cdot x,z]}\leq \beta_x)$. 
Let $A_x$ be $\{y\in A:(T_y,\beta_y)=(T_x,\beta_x)\}$. 
Note that for any such $y$, $\forall^*g\in G,\exists f \in L_{\beta_y}[g\cdot y,z]$ with 
$f$ witnessing $g\cdot y \in A$, by the absoluteness of illfoundedness for trees;
thus by the remarks following 1.3, 
$A_x$ is uniformly Borel in any real coding $(T_x,\beta_x)$. 
Since for any such $y$ and $G$, $G_{g\cdot y}\in L_{\beta_x}[g\cdot y,z]$ 
we can conclude that $E^X_G|_{A_x}$ is 
uniformly Borel in any code for $(T_x,\beta_x)$, by the same appeal to 1.8 made in the course of 3.1. 

Now if for some $x$ $E_0\sqsubseteq_c E^X_G|_{A_x}$ then we are done. 
Otherwise, as in the proof of 3.1, for each $x$ there is some $\theta_x$ reducing 
$E^X_G|_{A_x}$ to id$(2^{<\omega_1})$ with $\theta_x$ definable from 
$(T_x,\beta_x)$, and thus $G$-invariant. 
Now we may let 
$\theta(x)=(T_x,\beta_x,\theta_x(x))$. Since $T_x$ can be identified with a countable subset of 
$\omega_1^{<\omega}$, $\theta$ can be reorganized to give a reduction of $E_G^X|_A$ into 
id$(2^{<\omega_1})$.$\Box$\\

Some orbit equivalence relations allow HC-invariants, while others refuse them. While it may be 
interesting to explore new 
types of invariants for these more unruly equivalence relations, it turns out that in one direction 
the search is futile: If a Polish group action allows H$(\kappa)$ invariants for some ordinal 
$\kappa$, then it allows HC-invariants. 

Before the proof of this theorem we need one more basic fact from the theory of determinacy. 
Rephrased for the one context in view it reads as:\\

3.4. Theorem(Becker-Kechris). Assume AD$^{L(\R)}$. 
Let $T$ be the usual tree for $(\Sigma^2_1)^{L(\R)}$. Then for $x$ a real and 
$A\subset \ubf{\delta}^2_1$ with $A$ $(\Sigma^2_1(x))^{L(\R)}$ in the codes, 
$A\in L[T,x]$ and is uniformly definable from $x$ and $T$ over this model. (See \cite{beckerkechris2}.)\\

Recall in what follows that $\ubf{\delta}^2_1$ -- 
by definition the supremum of the $\Ubf{\Delta}^2_1$ 
prewellorders of $\R$ -- is equal to the least ordinal $\delta$ such that 
$L_{\delta}(\R)\prec_{\Ubf{\Sigma}^2_1} L(\R)$ . This standard fact follows as in the proof that $\delta^1_2$ 
is least $\delta$ such that $L_{\delta}\prec_{\Sigma_1} L$.\\

3.5. Theorem. Assume AD$^{L(\R)}$. Let $G$ and $X$ be as usual. 
Let $\kappa$ be an ordinal. 
If $E^X_G\leq_{L(\R)}$id(H$(\kappa)$), then $E^X_G\leq_{L(\R)}$id(H$(\omega_1)$). 

Proof. Choose $z$ a real coding the group and the action. 
If there exists a reduction into id(H$(\kappa)$), $\theta$, then the least 
$\alpha$ such that $L_{\alpha}(\R)$ satisfies DC+AD+ZF without Power Set and that such a reduction 
exists will be less than $\ubf{\delta}^2_1$. 
Without loss of generality, $\theta$ is definable from $z$ over $L_{\alpha}(\R)$. 
By leastness of $\alpha$ we have that $L_{\alpha}(\R)$ 
is the Skolem hull of its reals.  
Note that there is a relation $R\subset X\times \R$ in $L(\R)$ 
such that:

(i) $\forall x\in X \exists y\in \R((x,y)\in R)$; 

(ii) $(x,y)\in R$ implies that some set $B_{x,y}\subset\alpha$ coding 
$\theta(x)\in$H$(\alpha)$ 
is uniformly definable over $L_{\alpha}(\R)$ from 
$y$.

Let $T$ be the tree for $\Sigma^2_1(z_0)$, derived from the scale 
$(\varphi_n)_{n\in\omega}$. 
Then we find a parameter $z_0$ that codes $z$ 
with $\varphi_0(z_0)\geq\alpha$. 
Then the set of $\beta$  in 
$B_{x,y}$ is  uniformly $\Sigma^2_1(x,y,z_0)$ in the codes for any $(x,y)\in R$. 

Since $R$ is $\Sigma^2_1(z_0)$, for any real $x$ there is 
some $y\in L[x,z_0,T]$ with $(x,y)\in R$, and then by 3.4, $B_{x,y}$ is in $L[x,z_0,T]$, 
as is therefore $\theta(x)$, and both are uniformly definable in $L[x,z_0,T]$ from 
$z_0$, $x$ and $T$. 

Now, as in the proof of 3.1, we may find $M_x$ such that 
\[\forall^*g\in G({\rm HOD}^{L[z_0,g\cdot x, T]}_{(z_0,T,\theta(x))}=M_x),\]
and $x\mapsto M_x$ is $G$-invariant. Since $g\cdot x$ is generic over 
$M_x$ for the Vopenka algebra for (OD$(z, T, \theta(x))$)$^{L[z,g\cdot x, T]}$, 
there will be some ${\Bbb B}_x\in M_x$, $b_x\in{\Bbb B}_x$ with 
\[b_x\Vdash_{{\Bbb B}_x}\theta(x(\dot{G}))=\theta(x).\]
Note that isomorphism type of 
${\Bbb B}_x$ is canonical, as it equals 
 of the OD$(z, T, \theta(x))^{L[z,g\cdot x, T]}$ 
subsets of $X$ under inclusion. 
Note then that $b_x$ is also definable, as the union  of all such conditions in the algebra. 
Note that ${\Bbb B}_x$ has size at most $(2^{2^{\aleph_0}})^{M_x}$.  
$M_x$ can be wellordered in $L(\R)$ since $\theta(x)$ allows a wellorder;  
thus we have that $(2^{2^{\aleph_0}})^{M_x}<\omega_1$.  
Thus $({\Bbb B}_x,b_x)$ is a countable structure whose isomorphism type describes $[x]_G$. 

So let $\bar{\theta}(x)$ be its Scott sentence. $x\mapsto \bar{\theta}(x)$ is 
$G$-invariant and provides a complete invariant of $[x]_G$ since the above mentioned 
Vopenka algebra describes the equivalence class, in that if $\pi:({\Bbb B}_x,b_x)\cong ({\Bbb B}_y,b_y)$ 
then for $H_0\subset {\Bbb B}_x$ sufficiently generic, $H_1=\pi"[H_0]$ the push out of $H_0$, 
we have $\theta(x)=\theta((x(H_0))=\theta(x(H_1))=\theta(y)$, 
and so $\theta(x)=\theta(y)$ and $xE_G^X y$ as required.$\Box$

\newpage

{\bf $\S$4. Infinitary logic and group actions}\\

This section summarizes the main points in the development of infinitary logic for 
descriptive set theory and $\infty$-Borel codes. These remarks are along the 
lines of \cite{harrington} and \cite{harringtonshelah}, but with particular emphasis on the 
context of Polish group actions. 

The results here are technical and should be considered folklore. They will form the background for 
$\S$5.\\

4.1. Lemma. Let $X$ be a Polish space, ${\cal B}$ a basis for the topology, and $C\subset X$ a 
closed subset, Then ${\cal B}^C=\{O\cap C: O\in {\cal B}\}\cup{\cal B}$ is a basis of a new topology on 
$X$. 

Proof. Since $X$ with the new topology is homeomorphic to the disjoint union of 
$C$ and $X\setminus C$, both of which are shown in \cite{kechris2} to be Polish in the 
relative topology.$\Box$\\

4.2. Definition. Let $X$ be a Polish space and ${\cal B}$ a basis. Let ${\cal L}({\cal B})$ be the 
propositional language formed from atomic propositions of the form $'\dot{x}\in U'$, for $U\in{\cal B}$. 
Let ${\cal L}_{\omega_10}({\cal B})$ be the infinitary version, obtained by closing under negation and 
countable disjunction and conjunction, and let ${\cal L}_{\infty 0}({\cal B})$ be the obtained by 
closing under arbitrary Boolean operations. 
$F\subset$ ${\cal L}_{\infty 0}({\cal B})$ is a {\it fragment} if it is closed under subformulas and 
the finitary Boolean operations of negation and finite disjunction and finite conjunction. 
For $\varphi\in$${\cal L}_{\infty 0}({\cal B})$, $F(\varphi)$, the {\it fragment generated by} 
$\varphi$, is the smallest collection fragment containing $\varphi$. 

For a point $x\in X$ and $\varphi\in$${\cal L}_{\infty 0}({\cal B})$, we can then define 
$x\models \varphi$ by induction in the usual fashion: If $\varphi='\dot{x}\in U'$ then $x\models \varphi$ 
if and only if $x\in U$; for $\varphi=\neg\psi$, $x\models \varphi$ if and only if it is not the case 
that $x\models \psi$; for $\varphi=\bigwedge\{\psi_i:i\in\Lambda\}$, 
$x\models\varphi$ if and only if for every 
$i\in\Lambda$ we have 
$x\models \psi_i$. 

For $F\subset{\cal L}_{\infty 0}({\cal B})$ a countable set closed under subformulas, 
we let 
$\tau(F)$ be the topology generated by ${\cal B}$ and 
all sets of the form $\{x\in X:x\models \varphi\}$, as $\varphi$ ranges over 
$F$.\\

4.3. Lemma. For $F\subset {\cal L}_{\omega_10}({\cal B})$ a countable fragment, 
$\tau(F)$ forms the basis of a Polish topology on $X$. 

Proof: This follows by 4.1 and induction on the complexity of the infinitary sentences in $F$. 
If $\psi=\neg\phi$ it follows by inductive assumption and 4.1. For $\psi=\bigvee_{\Lambda}
\psi_i$ it is trivial since we are simply adding a new open set to the basis. At limit stages of 
the construction we may use that increasing countable unions of Polish topologies are again 
Polish -- a classical fact that is recalled in \cite{sami} and \cite{kechris2}.$\Box$\\

Note then that if $F$ is a fragment of ${\cal L}_{\infty 0}({\cal B})$ 
then in any generic extension in which $F$ becomes countable it must generate 
a Polish topology. We will frequently have cause to consider Polish spaces and continuous 
Polish groups both in $\V$ and through future generic extensions. 
This is reasonable, 
since all the relevant statements of the 
form '$X$ is a Polish $G$-space' are $\Ubf{\Pi}^1_1$, and hence absolute. 

The next lemma merely makes the point that we may find the Vaught transform of a $\varphi\in 
{\cal L}_{\infty 0}({\cal B})$ in a manner that is effective. The proof is a rephrasing of 
the usual proof that the Vaught transform of a Borel set is again Borel.\\

4.4. Lemma. Let $G$ be a Polish group, $X$ a Polish $G$-space, ${\cal B}$ a countable 
basis for $X$, ${\cal B}_0$ a countable basis for $G$. Then to each $\varphi\in {\cal L}_{\infty 0}({\cal B})$ 
and $V\in {\cal B}_0$ we may assign 
a formula $\varphi ^{\Delta V}\in {\cal L}_{\infty 0}({\cal B})$ such 
that:

(i) $(V, \varphi)\mapsto \varphi ^{\Delta V}$ is uniformly 
$\Ubf{\Delta}_1$ in any 
parameter coding $X$, $G$, the action, and the bases;

(ii) the fragment generated by 
$\varphi ^{\Delta V}$ has the same cardinality as 
the fragment generated by $\varphi$, and in fact they have approximately the same 
logical complexity; 

(iii) in all generic extensions $\V[H]$ of $\V$ in which $\varphi\in ({\cal L}_{\omega_10}({\cal B}))^{\V[H]}$ 
we have 
\[\{x\in X: x\models \varphi ^{\Delta V}\}=
\{x\in X:\exists^* g\in V(g\cdot x\models \varphi)\}.\]

Proof. By the usual type of induction on the logical complexity of $\varphi$.$\Box$\\

Note here that the calculation of whether $x\models \varphi$ is absolute to any model containing $x$ and 
$\varphi$. The statement of 4.4 gives that the assignment $(V, \varphi)\mapsto \varphi ^{\Delta V}$ 
will be $\Ubf{\Delta}^1_2$ when restricted to $\varphi\in {\cal L}_{\omega_10}$, 
since $\Ubf{\Delta}_1^{HC}=\Ubf{\Delta}^1_2$.\\

4.5. Lemma. Let $G$, $X$, ${\cal B}$, 
${\cal B}_0$ be as in 4.4. 
Then to each $\varphi\in {\cal L}_{\infty 0}({\cal B})$ 
and $g\in G$ we may assign 
a formula $\varphi ^{g}\in {\cal L}_{\infty 0}({\cal B})$ such 
that:

(i) $(g, \varphi)\mapsto \varphi ^{g}$ is uniformly 
$\Ubf{\Delta}_1$ in any 
parameter coding $X$, $G$, the action, and the bases;

(ii) the fragment generated by 
$\varphi ^{g}$ has the same cardinality as 
the fragment generated by $\varphi$; 

(iii) in all generic extensions $\V[H]$ of $\V$ in which $\varphi\in ({\cal L}_{\omega_10}({\cal B}))^{\V[H]}$ 
we have 
\[\{x\in X: x\models \varphi ^{g}\}=
\{x\in X:g\cdot x\models \varphi\}.\]

Proof. Again by induction on $\varphi$.$\Box$\\

The next lemma can be contrasted with the notion of $\infty$-Borel from $\S 2$. 
In effect the lemma states that every $\infty$-Borel code is representible by some infinitary 
$\varphi\in {\cal L}_{\kappa 0}({\cal B})$, for some ordinal $\kappa$.\\

4.6. Lemma. Let $X$ be a Polish space and ${\cal B}$ be a basis. Then for $\alpha$ an ordinal, 
$S\subset \alpha$, and $\psi\in{\cal L}(\in)$ be a formula in set theory, there is a corresponding 
$\varphi(\alpha, S,\psi)\in {\cal L}_{\infty 0}({\cal B})$ such that 

(i) in all generic extensions, $\{x\in X:x\models \varphi(\alpha, S,\psi)\}=\{x\in X: 
L_{\alpha}[S,x]\models \psi(x, S)\}$; 

(ii) the transitive closure of $\varphi(\alpha, S,\psi)$ has cardinality $|\alpha|+\aleph_0$, 
so that for $G\subset$Coll$(\omega,\alpha)$ $\V$-generic, $\varphi(\alpha, S,\psi)\in 
{\cal L}_{\omega_1 0}({\cal B})$; 

(iii) the assignment $\alpha, S,\psi\mapsto \varphi(\alpha, S,\psi)$ is $\Delta_1$ in any parameter coding the 
space and the basis. 

Proof. For the purposes of this argument, let us suppose that $X=2^{\omega}$ and that 
${\cal B}$ is the usual basis, $\{\{x\in X:x\supset s\}s\in 2^{<\omega}\}$. 
The more general case can be handled similarly, with the details depending on our 
precise manner of coding and constructing from points in an arbitrary Polish space. 

The proof is by induction on the complexity of the $(\alpha, S,\psi)$, with the 
base case corresponding to $\alpha=\omega$, $S\subset\omega$, $\psi\in\Sigma_1$. 
Then the set $\{x\in X: 
L_{\alpha}[S,x]\models \psi(x, S)\}$ is open, and trivially representable by a 
formula in ${\cal L}_{\infty 0}({\cal B})$. 
Carrying the induction through complementation is immediate. 
If $\psi\in\Sigma_{n+1}$ and we have proved the lemma for $\Pi_n$ over 
$L_{\alpha}[S,x]$, then 
\[\{x\in X:L_{\alpha}[S,x]\models \psi(S,x)\}=\bigcup\{\{x\in X:L_{\alpha}[S,x]\models
 \overline{\psi}(S,x,\beta)\}:\beta\in\alpha\}\] 
for some $\overline{\psi}\in\Pi_n$, and thus follows by inductive step. 

The limit case of the induction corresponds to considering 
$\Sigma_1$ over $L_{\alpha}[S,x]$ for $\alpha>\omega$, and follows as 
in the successor step for $\Sigma_{n+1}$ above.$\Box$\\

4.7. Theorem(Becker-Kechris). Let $G$ be a Polish group and $X$ a Polish $G$-space, ${\cal B}$ a 
basis for $X$, ${\cal B}_0$ a basis for $G$, $G_0\subset G$ be a countable 
dense subgroup. Let ${\cal C}$ be a collection of Borel sets in $X$ such 
that

(i) ${\cal C}$ is an algebra -- in other words, closed under finite Boolean operations; 

(ii) ${\cal C}$ is closed under translation by elements in $G_0$; 

(iii) ${\cal C}$ is closed under Vaught transforms from ${\cal B}_0$, so that for 
$C\in {\cal C}$ and $U\in {\cal B}_0$, $C^{* U}, C^{\Delta U}\in{\cal C}$; 

(iv) ${\cal C}$ forms the basis of a Polish topology on $X$. 

Then: $\{C^{\Delta U}:C\in{\cal C}, U\in{\cal B}_0\}$ forms the basis of a Polish 
topology on $X$ under which it remains a Polish $G$-space. (See \cite{beckerkechris}.)\\

Note that 
(ii) guarentees that for any $C_1..C_l\in{\cal C}$, 
$U_1,...U_l\in{\cal B}_0$, $C_1^{\Delta U_1}\cap...\cap C_l^{\Delta U_l}$ is the 
union of sets in $\{C^{\Delta U}:C\in{\cal C}, U\in{\cal B}_0\}$.\\

4.8. Lemma. Let $X$, ${\cal B}$, be as 4.7, let $E$ be a $\Ubf{\Sigma}^1_2$ equivalence 
relation on $X$, and 
let 
${\Bbb P}$ be a forcing notion, $p\in{\Bbb P}$, and $\sigma$ a term for  
$\V^{\Bbb P}$ a term for an element in $X$, such that 
\[(p,p)\Vdash_{{\Bbb P}\times{\Bbb P}}\sigma[\dot{G}_l]E\sigma[\dot{G}_r],\]
where, as usual, $\dot{G}_l$ and $\dot{G}_r$ refer to the generic objects on the 
left and right copies of ${\Bbb P}$. 

Then there is a set $F\subset$${\cal L}_{\infty 0}({\cal B})$ and a $\varphi\in F$ 
such that:

(i) $\tau(F)$ generates a Polish topology on $X$ in any generic extension in which $F$ 
becomes countable; 

(ii) for any generic $H\subset{\Bbb P}$, $\sigma[H]\models \varphi$; 

(iii) for any $x_1, x_2 $ in any generic extension of $\V$ with $x_1, x_2\models \varphi$ 
we have $x_1Ex_2$. 

Proof. Following 4.6, we can certainly find $\varphi$ such that for any $x\models \varphi$ one must 
have that $x$ is $\V$-generic forcing for the factor forcing that introduces $\sigma[H]$, 
below $p\in{\Bbb P}$, 
for some $\V$-generic $H\subset {\Bbb P}$, with $H$ possibly only appearing in a further 
generic extension of $V[x]$. By closing $\varphi$ under subformulas and finite Boolean 
operations, we obtain a Polish topology by 4.3. Thus we have (i) and (ii). 

So now suppose that $x_1, x_2\models \varphi$. Then we can generically find $H_1$, $H_2$ 
that are $\V$-generic below $p$ with $\sigma[H_1]=x_1$, $\sigma[H_2]=x_2$. Then by choosing 
$H_3\subset {\Bbb P}$ sufficiently generic below $p$, and setting $x_3=\sigma[H_3]$ we 
get that $H_1\times H_3$ and $H_1\times H_2$ are both $\V$-generic. 

Then by the assumptions on ${\Bbb P}$, $p$, and $\sigma$, 
\[x_1Ex_3,\]
\[x_2E x_3,\]
and thus
\[x_1Ex_2.\]$\Box$\\
 
In the context of Polish group actions  
4.7 suggest a refinement.\\

4.9. Corollary. Let $G$, $X$, ${\cal B}$, $G_0$ be as 4.7, and let 
${\Bbb P}$  $p\in{\Bbb P}$, and $\sigma$ be as in 4.8. 

Then there is a set $F_0\subset$${\cal L}_{\infty 0}({\cal B})$ and a $\varphi_0\in F_0$ 
such that:

(i) $\tau(F_0
)$ generates a Polish topology on $X$ in any generic extension 
$\V[H]$ in which $F_0$ becomes countable, 
{\it and} $(X, \tau(F_0))$ remains a Polish $G$-space; 

(ii) for any generic $H\subset{\Bbb P}$ 
\[\forall g\in G(g\cdot 
\sigma[H]\models \varphi_0);\] 

(iii) for any $x_1, x_2 $ in any generic extension of $\V$ with $x_1, x_2\models \varphi_0$ 
we have $x_1E^X_Gx_2$. 

Proof. First let ${\Bbb P}^*$ be the forcing notion of ${\Bbb P}$ followed by the version 
of 
Cohen forcing obtained by 
using the basic open sets in $G$ to create a generic group element. 
The let $\tau$ be the term in ${\Bbb P}^*$ for $\dot{g}\cdot \sigma[\dot{H}]$, 
where $\dot{g}$ names the generic group element and $\dot{H}$ denotes the 
generic on ${\Bbb P}$, 
and let $q=\langle p,1\rangle$ be the condition in ${\Bbb P}^*$ obtained by 
insisting that $p$ be in the generic $\dot{H}$. 

Then ${\Bbb P}^*$, $q\in {\Bbb P}^*$, and $\tau$ continue to satisfy the assumptions 
of 4.8, but we have engineered the futher result that if $(h,H)$ is a generic on 
${\Bbb P}^*$, and $x=\tau[(H,h)]$, then in any in any future generic extension in 
which $(2^{{{\Bbb P}^*}})^{\V}$ becomes countable we have that 
$\forall^*g\in G((gh, H)$  is ${\Bbb V}$ generic for $
{\Bbb P}^*$ below $\langle p,1\rangle $).

Now we choose $F$ and $\varphi$ 
as in 4.8, for ${\Bbb P}^*$ and $\tau$, but taking enough care to ensure closure under 
$G_0$ translation and $\Delta$-Vaught transforms with respect to ${\cal B}_0$. 
This can certainly be achieved by 4.4 and 4.5. 
So then we obtain (i), (ii), and (iii) as in 4.8, but with the further condition that 
in any generic extension of $\V$ containing $x$ in which $(2^{{\Bbb P}^*})^{\V}$ becomes countable
\[x\models \varphi\Rightarrow \exists^*g\in G(g\cdot x\models \varphi),\] 
and so in the notation of 4.4, 
\[x\models \varphi\Rightarrow x\models \varphi ^{\Delta G}.\]

Claim: In all generic extensions, $\varphi^{\Delta G}$ is $G$-invariant. 

Note that any generic extension in which 
$\varphi$ is in ${\cal L}_{\omega_1 0}({\cal B})$ we will have that $x\models \varphi^{\Delta G}$ if and 
only if there is a non-meager collection of group elements $g\in G$ such that 
$g\cdot x\models \varphi$. Thus $x\models \varphi^{\Delta G}$ if and only if 
$g_0\cdot x \models \varphi^{\Delta G}$ for any $g_0\in G$, in this model, and hence also 
in $\V[x,g_0]$ -- since the calculation of $x\models \varphi$, $g_0\models \varphi^{\Delta G}$ are 
absolute to $\V[x,g_0]$.($\Box$Claim)

So if we now follow 4.7 and let 
$F_0$ be $\{\psi^{\Delta U}: \psi \in F, U\in {\cal B}_0\}$, then we can take 
$\varphi_0=\varphi^{\Delta G}\in F_0$. 
Then this is as required since $\varphi^{\Delta G}$ is $G$-invariant.$\Box$\\

\newpage
{\bf $\S$5. Becker's theorem}\\

5.1. Lemma. Let $G$ be a Polish group admitting a left complete 
invariant metric, and $X$ a Polish $G$-space, 
and ${\cal B}$ a countable basis for $X$. 
Let $x\in X$. 
Let ${\cal M}$ be a class inner model of ZF+DC, with $X$, $G$, and the action 
existing in ${\cal M}$, in the sense of being coded by a parameter in 
${\cal M}$. 
Let $F\subset {\cal L}_{\infty 0}({\cal B})$ be a set 
closed under subformulas such that 

(i) $F\in {\cal M}$; 

(ii) in some generic extension of ${\Bbb V}$ $\tau(F)$ generates a Polish topology on 
$X$, including the orginal topology, with $(X, \tau(F))$ a Polish $G$-space; 

(iii) $[x]_G$ is $\tau(F)$-open. 

Then: $[x]_G\bigcap {\cal M}\neq\emptyset.$

Proof. Note that if $H\subset$Coll$(\omega, F)$  is ${\cal M}$-generic, then 
 ${\cal M}[H]\bigcap [x]_G\neq\emptyset$, and so  
$\{\varphi\in F:\exists x_0\in[x]_G(x_0\models \varphi)\}\in{\cal M}$; the task is to use 
this kind of information to find a representative of the orbit. 
Note also that the natural map from $G$ to $[x]_G$, $g\mapsto g\cdot x$,  
is open and continuous by 1.9. 

Let ${\cal B}_0$ be a countable basis for $G$ in ${\cal M}$. 
Let $d_r$ be a right invariant complete metric on $G$.  
Let $d$ be a complete metric on $X$. 
Then for any $V\in{\cal B}_0$ with $1_G\in V$ 
we let 
\[{\cal B}(V)=\{\varphi \in F: \forall x_0, x_1\in[x]_G
(x_0, x_1\models \varphi\Rightarrow \exists g\in V(g\cdot x_0=x_1))\}.\]
By 1.9, this is non-trivial for $V\neq\emptyset$; 
by the genericity of $[x]_G$, the function $V\mapsto {\cal B}(V)$ exists in 
${\cal M}$. 

Thus using DC in ${\cal M}$, we may find sequences 

\[(V_n)_{n\in{\Bbb N}}\subset {\cal B}_0,\]
\[(U_n)_{n\in{\Bbb N}}\subset {\cal B},\]
\[(\varphi_n)_{n\in{\Bbb N}}\subset F,\]
and 
\[(z_n)_{n\in{\Bbb N}}\subset X\]
such that if $H\subset$Coll$(\omega, F)$ is ${\cal M}$-generic then 
${\cal M}[H]$ satisfies 
\[d_r(V_n)<2^{-n},\]
\[d(U_n)<2^{-n},\]
\[\exists x_0
\in [x]_G\cap U_{n+1}
(x_0\models\varphi_n),\]
\[z_n\in U_n,\]
\[\varphi_n\in {\cal B}(V_n),\]
\[y\in X(y\models \varphi_{n+1}\Rightarrow y\in U_{n+1}\wedge y\models \varphi_n),\]
\[\overline{U}_{n+1}\subset U_n.\]

The above assignments exist already in ${\cal M}$ since 
it has access to the function $V\mapsto {\cal B}(V)$. 
Since $(z_n)_{n\in{\Bbb N}}$ is Cauchy in ${\cal M}$, we can find 
$z_{\infty}\in{\cal M}$ such that $z_n\rightarrow z_{\infty}$ as $n\rightarrow \infty$. 

Meanwhile in ${\cal M}[H]$ choose $(x_n)_{n\in{\Bbb N}}$ such that $x_n\in[x]_G\cap U_n$ and 
$x_n\models \varphi_n$. By the definition of ${\cal B}(V_n)$ there are group elements 
$g_n\in V_n$ 
such that 
\[g_n \cdot x_n=x_{n+1}.\]

Note that $d_r(g_n, 1_G)<2^{-n}$ implies $d_r(g_n\cdot g_{n-1}\cdot...g_0, g_{n-1}\cdot...g_0)<2^{-n}$. 
Thus if we set $h_n=g_n\cdot g_{n-1}\cdot...\cdot g_0$ then 
$(h_n)_{n\in{\Bbb N}}$ is a Cauchy sequence in $G$ with respect to 
$d_r$. So there is some $h_{\infty}$ that is the limit $(h_n)_{n\in{\Bbb n}}$, and 
$x_{\infty}=h_{\infty}\cdot x_0$lim$_{\Bbb N}x_n$. 
Since $d(x_n, z_n)<2^{-n}$ we get $x_{\infty}=z_{\infty}\in{\cal M}\cap [x]_G$, as required.$\Box$\\

5.2. Corollary. Let ${\cal M}$ be a class inner model of ZF+DC. 
Let $X$ be a Polish $G$-space, with both objects along with the action existing in ${\cal M}$. 
Let $({\Bbb P},\sigma, p)\in{\cal M}$ 
such that $\sigma$ is term for the ${\cal M}^{\Bbb P}$ such that 
\[(p,p)\Vdash_{{\Bbb P}\times{\Bbb P}} \sigma[\dot{G}_l] E^X_G \sigma[\dot{G}_r].\]

Then there is some $y\in{\cal M}$ such that 
\[p\Vdash_{\Bbb P}\sigma[\dot{G}] E^X_G y.\]

Proof. We may as well assume that ${\Bbb V}$ has a representative of the generic equivalence class, 
since otherwise we may replace ${\Bbb V}$ by ${\Bbb V}[H]$ for some suitably generic $H$. 
Then the theorem follows by 5.1 and 4.9.$\Box$\\

It follows then from the results of \cite{gao} that 5.2 characterizes when a closed subgroup 
of $S_{\infty}$ is cli, in 
that if $G$ is a closed subgroup of the symmetric group that does not 
admit a left invariant complete metric then there is a Polish $G$-space $X$ and 
$\sigma$ a term for the forcing notion ${\Bbb P}=$Coll$(\omega, \omega_1)$ such 
\[{\Bbb P}\times{\Bbb P}\Vdash \sigma[\dot{G}_l]E^X_G \sigma[\dot{G}_r],\]
and for all $x\in X$
\[{\Bbb P}\Vdash \neg(\sigma[\dot{G}]E^X_G).\]
I do not know if this characterization succeeds for arbitrary Polish groups.\\

5.3. Theorem. TVC(cli,$\Ubf{\Sigma}^1_1$) -- which is to say, if $G$ is a cli Polish group, 
$X$ a Polish $G$-space, $A\subset X$ $\Ubf{\Sigma}^1_1$, 
then either $|A/G|\leq\aleph_0$ or there is a perfect set $P\subset A$ such that any two 
elements in $P$ are $E^X_G$-inequivalent. 

Proof. If $A$ has uncountably many orbits, then as shown in \cite{sami}, this is 
$\Ubf{\Pi}^1_2$: For all $x_0, x_1,...\in X$ and $F_0, F_1,...$ closed subsets of 
$G$, either: 

(i) there is some $n$ such that $F_n$ is not the stabilizer of $x_n$; or, 

(ii) using $F_n$ to obtain a  uniform calculation of $[x_n]_G$ as a Borel 
set, there is some $x\in A\setminus\bigcup_n[x_n]_G$. 

Thus through all generic extensions there will always be uncountably many 
orbits in $A$. Thus for ${\Bbb P}=$Coll$(\omega,\kappa)$, $\kappa$ sufficiently big, 
there will some term $\sigma$ with 
\[{\Bbb P}\Vdash \forall x\in {\Bbb V}(x\not\in[\sigma[\dot{G}]]_G).\]
Note also that the 
existence of a perfect set of inequivalent reals in $A$ is $\Ubf{\Sigma}^1_2$ 
since, since granted perfect $P\subset A$, we can find perfect $P_0\subset P$ 
and continuous $f$ with domain $P_0$ such that for any $x\in P_0$, 
$f(x)$ witnesses $x\in A$. So again the absence of a perfect set of orbits in 
$A$ will again be absolute.

The usual sort of diagonalization arguments -- as can be found in \cite{ha1}, \cite{harringtonshelah}, or 
\cite{stern} -- show that either there is a term a condition $p\in{\Bbb P}$ that decides the 
equivalence class, in the sense that 
\[(p,p)\Vdash_{{\Bbb P}\times{\Bbb P}}\sigma[\dot{G_l}]E^X_G \sigma[\dot{G_r}],\] 
or there is a generic extension with a perfect set of inequivalent reals: In a 
generic extension choose $(D_n)_{n\in{\Bbb N}}$ an enumeration of the open dense subsets 
of ${\Bbb P}\times{\Bbb P}$, and then choose $(p_s)_{s\in 2^{<{\Bbb N}}}$ so that for $s\neq t\in 
2^n$, $(p_s, p_t)\in D_n$, and $(s0, s1)\Vdash_{{\Bbb P}\times{\Bbb P}}
\neg(\sigma[\dot{G_l}]E^X_G \sigma[\dot{G_r}])$; then for $G_x$ being the filter generated by 
$\{p_{x|n}:n\in{\Bbb N}\}$, any two distinct $x,y\in 2^{{\Bbb N}}$ give rise to 
$\sigma[\dot{G_x}]$ and $\sigma[\dot{G_y}]$ that are inequivalent; performing this with enough 
care we actually do finish with a perfect set of orbit inequivalent points in $A$. 
Then by absoluteness we obtain this in the ground model ${\Bbb V}$ 
and we are finished. 

So suppose instead that there is a condition $p$ deciding the equivalence class. 
For $H\subset {\Bbb P}$ ${\Bbb V}$-generic, $x=\sigma[H]$, 
in ${\Bbb V}[H]$ we can apply 5.2 to ${\cal M}={\Bbb V}$, 
and obtain that $[x]_G$ has a representative in ${\Bbb V}$, 
contradicting the 
assumption on ${\Bbb P}$ and $\sigma$.$\Box$\\

5.4. Theorem. Let $G$ be a Polish group with a left invariant complete metric 
acting continuously on a Polish space $X$, and let $A\subset X$ be $\Ubf{\Sigma}^1_1$. Then either 

(I) there is a $\Ubf{\Delta}^1_2$ 
$\theta:A\rightarrow 2^{\omega}$ such that for all $x_1,x_2\in A$ 
\[\exists g\in G(g\cdot x_1=x_2)\Leftrightarrow \theta(x_1)=\theta(x_2)\]
or 

(II) there is a Borel $\theta: {\Bbb R}\rightarrow A$ such that for all $r_1, r_2\in {\Bbb R}$ 
\[r_1-r_2 \in {\Bbb Q} \Leftrightarrow \exists g\in G(g\cdot \theta(r_1)=\theta(r_2)).\]

Proof. Let $A=\{x\in X:\exists y\in \omega^{\omega}B(x,y)\}$, for some Polish space $B\subset 
X\times \omega^{\omega}$. 
Define $E$ on $B$ by $(x_1,y_1)E(x_2, y_2)$ if and only if 
$x_1E^X_Gx_2$. This is a $\Ubf{\Sigma}^1_1$ equivalence relation such that through all generic 
extensions every equivalence class is Borel. 

Now we follow \cite{hjorthkechris}. One case is that $E_0\sqsubseteq_c E$, when we are quickly 
finished. Alternatively, we obtain a $\Ubf{\Delta}^1_2$ in the codes reduction into $2^{<\omega_1}$, 
call it $\theta$. Then it is $\Ubf{\Pi}^1_2$ to assert that 
\[\forall x_1, x_2\in X(\theta(x_1)=\theta(x_2)\Rightarrow x_1 E^X_G x_2.\,\]
and thus absolute. Let $z$ be a real coding the action, and any parameters used in the definition 
of $\theta$. 

Then for all $(x,y)\in B$ there must be a representative of the $E$-equivalence class of 
$(x,y)$ in any generic extension of $L[\theta(x,y), z]$ in which $\theta(x,y)$ is countable, 
by absoluteness of $\Ubf{\Sigma}^1_2$. Thus $[x]_G$ will be generic over 
$L[\theta(x,y),z]$, and thus by 5.2 there will be some $x_0\in[x]_G\bigcap 
L[\theta(x,y), z]$. 

So now we can define $\theta_0
:A\rightarrow X$ be letting $\theta_0(x)$ be the first real under the canonical wellorder in 
$L[\theta(x,y),z]$ with $xE^X_G x_0$. This gives a reduction of $E^X_G|_A$ to id$(X)$ which 
can in turn be reorganized to give a reduction into id$(2^{\omega})$ or id$(\R)$.$\Box$\\

Combining these ideas with the methods of $\S 3$ it can be shown that:\\

5.5. Theorem(AD$^{L(\R)}$). Let $G$ be a cli Polish group and let
$X$ be a Polish 
$G$-space. Let $A\subset X$ be in ${L(\R)}$. 
Then either: 

(I) $E^X_G|_A\leq$id$(2^{\omega})$; or, 

(II) $E_0\sqsubseteq_c E^X_G|_A$.\\

The observation that underlies the proof of 5.5 is that if we are in case (I) of 3.2, as 
witnessed by $\theta:A\rightarrow 2^{<\omega_1}$ in $L(\R)$, and if $S\subset\alpha$ is an $\infty$-Borel 
code for $\{(x,w_0,w_1 ,i): x\in X, w_0, w_1$ code $\alpha<\delta, \theta(x)\in 2^{\delta}, \theta(x)(\alpha)
=i\}$, then $\theta(x)\in L[S, x]$; in the notation of the proof of 3.1, $[x]_G$ will be generic 
over $M^x_S[\theta(x)]$; thus $[x]_G\bigcap M^x_S[\theta(x)]\neq \emptyset$ by 5.1. 

The next theorem states that the orbits of a cli group cannot be used to code countable sets of 
reals. Since this is one of the simplest equivalence relations induced by the symmetric group, 
and in some ways appears distinctive of this group, 
the result underscores the divergence between cli group actions and arbitrary orbit equivalence relations 
induced by $S_{\infty}$.\\

5.6. Theorem. Let $Y=\R^{\Bbb N}$ and 
define $E$ by $(y_n)_{n\in{\Bbb N}}E(x_n)_{n\in{\Bbb N}}$ if and only $\{y_n:n\in{\Bbb N}\}=
\{x_n:n\in{\Bbb N}\}$ -- so this is the orbit equivalence relation induced by the $S_{\infty}$-action 
of $(g\cdot \vec x)(n)=x(g^{-1}(n))$ for $\vec x \in Y$ and $g\in S_{\infty}$. 
Then there is no cli Polish group $G$ Polish $G$-space $X$ with $E\leq_B E^X_G$. 

Proof. Instead suppose $\theta:Y\rightarrow X$ performs a Borel reduction. 
Note that this statement is $\Ubf{\Pi}^1_2$: 
\[\forall y_1, y_2\in Y(y_1E y_2\Leftrightarrow \exists g\in G(g\cdot\theta(y_1)=\theta(y_2));\]
and hence it would be absolute through all generic extensions. 
 
Let ${\Bbb P}$ be the forcing to collapse $2^{\aleph_0}$ to $\omega$, and let $\sigma[\dot{G}]$ 
denote the term in ${\Bbb V}^{\Bbb P}$ an element of $Y$ that enumerates every real once. 
Thus 
\[{{\Bbb P}\times{\Bbb P}}\Vdash \sigma[\dot{G}_l] E \sigma[\dot{G}_r].\]
Thus if we let $\sigma_0[\dot{G}]$ be the term for $\theta(\sigma[\dot{G}])$. 
By the absoluteness of the assumptions on $\theta$ 
\[({{\Bbb P}\times{\Bbb P}}\Vdash \sigma_0[\dot{G}_l] E^X_G \sigma_0[\dot{G}_r]).\]
Thus by 5.2 there is some $x\in X$ with ${\Bbb P}\Vdash \sigma_0[\dot{G}] E^X_G x$. 
However, since ${\Bbb V}^{\Bbb P}\models \exists y\in Y(\theta(y)E^X_G x)$, this must 
hold in already in ${\Bbb V}$ by the absoluteness of $\Ubf{\Sigma}^1_2$. So fix 
$y\in Y$ with $\theta(y)E^X_G x$. Then again by the absoluteness of the assumptions on 
$\theta$, ${\Bbb P}\Vdash \sigma[\dot{G}] E y$. 

This is absurd, since any such $y$ would need to enumerate ${\Bbb R}$ in order type $\omega$.$\Box$\\

While a similar result is proved for abelian groups in \cite{hjorth2} without an appeal to 
metamathematics, the only known proof of 5.6 uses forcing.\\
\newpage

{\bf $\S$6 Knight's model}\\

6.1. Definition. Let $\sigma\in{\cal L}_{\omega_1\omega}$, for ${\cal L}$ some countable language, 
which we may assume without loss of generality to be relational. 
Then Mod$(\sigma)$ is the set of all models of $\sigma$ whose underlying set is ${\Bbb N}$. 
We let $\tau(\sigma)$ be the topolgy generated by sets of the form $\{M\in$Mod$(\sigma):M\models 
\psi(n_0, n_1, ...,n_k)\}$ where $(n_0,...n_k)$ a finite sequence of natural numbers 
and $\psi$ is a formula in the fragment generated by 
$\sigma$, in the sense that it is in the smallest collection of formulas containing $\sigma$ and 
closed under subformulas, substitutions, and the first order operations of negation, finite disjunction, finite 
conjunction, and existential quantifiers. 
Mod${\cal L}$ is the collection of ${\cal L}$ models on ${\Bbb N}$ with the topology generated by first 
order logic. 

We then let $S_{\infty}$ act on Mod$(\sigma)$ by 
\[(g\cdot M)\models R(n_0,...,n_k) \Leftrightarrow M\models R(g^{-1}(n_0),...,g^{-1}(n_k)),\]
for any $R\in{\cal L}$, $(n_0,...,n_k)$ a finite sequence in ${\Bbb N}$. 
The equivalence relation $E_{S_{\infty}}$ induced by this action on Mod$(\sigma)$ is frequently 
denoted by $\cong|$Mod$(\sigma)$.\\

6.2. Lemma(folklore). For any $\sigma\in{\cal L}_{\omega_1\omega}$, (Mod$(\sigma), \tau(\sigma)$) 
is a Polish $S_{\infty}$-space.\\

6.3. Definition. For $M$ a model and $\vec a\in M^{<\omega}$ one defines the canonical $\alpha$ type 
of $\vec a$, $\varphi_{\alpha}^{\vec a, M}\in{\cal L}_{\infty\omega}$, 
by induction on $\alpha$: $\varphi_0^{\vec a, M}$ is the infinitary formula expressing the 
quantifier free type of $\vec a$ in $M$. 
\[\varphi_{\alpha+1}^{\vec a, M}=\varphi_{\alpha}^{\vec a, M}\bigwedge _{b\in M}\exists x
\varphi_{\alpha}^{{\vec a}b, M}\wedge\forall x\bigvee_{b\in M}
\varphi_{\alpha}^{{\vec a}b, M}.\]
At limit stages we take intersections. 

The {\it Scott height} of $M$ is the least $\gamma$ such that for all $\vec a$, 
$\varphi_{\delta}^{\vec a, M}$ determines $\varphi_{\delta +1}^{\vec a, M}$. 
The {\it Scott sentence} of $M$, $\varphi_M\in {\cal L}_{\infty\omega}$ states what 
$\gamma$-types exist for $\gamma$ the Scott height and that this {\it is} the Scott height. 
As in \cite{keisler}, two countable models are isomorphic if and only if they have the 
same Scott sentence.\\
 
6.4. Definition. Let $M$ be a countable model with underlying set ${\Bbb N}$. Then Aut$(M)=
\{g\in S_{\infty}:g\cdot M=M\}$.\\

6.5. Theorem(folklore). Let $G$ be a subgroup of $S_{\infty}$. Then $G$ is closed in $S_{\infty}$ if and 
only if $G=$Aut$(M)$ for some countable $M$ with underlying set ${\Bbb N}$.\\

The authors of \cite{beckerkechris} noticed that this allows a curious analogue in the context of 
Polish group actions.\\

6.6. Theorem(Becker-Kechris). Let $G=$Aut$(M)$ be a closed subgroup of $S_{\infty}$; let 
${\cal L}$ be the language of $M$. Let $X$ 
be a Polish $G$-space. Then there is a language 
${\cal L}'\supset{\cal L}$ extending the language of $M$ 
and $\sigma\in {\cal L}'$ such that $\sigma\Rightarrow \varphi_M$ and $|X/G|$ is 
Borel equivalent to Mod$(\sigma)$, in the sense that there are $\theta:X\rightarrow$Mod$(\sigma)$ 
and $\rho:$Mod$(\sigma)\rightarrow X$ such that:

(i) $\theta$ witnesses $E^X_G\leq_B \cong|$Mod$(\sigma)$; 

(ii) $\rho$ witnesses $\cong|$Mod$(\sigma)\leq_B E^X_G$; and 

(iii) these are orbit inverses to one another in the sense 
that for all $x\in X$, $xE^X_G(\rho\circ\theta(x))$. 

Proof(sketch). Let $({\cal O}_m)_{m\in{\Bbb N}}$ be a countable basis for $X$. 
We may associate to each $x\in X$ the model $M_x$, with relations $(R_{m,k})_{m,k\in{\Bbb N}}$, 
where for $(n_1,...,n_k)$ a $k$-tuple in ${\Bbb N}$, 
\[M_x\models R_{m,k}(n_1,...,n_k) \Leftrightarrow \forall^*g\in G(g(n_1)=0\wedge...\wedge g(n_k)=k-1\Rightarrow 
g\cdot x\in {\cal O}_m).\] 
It is shown in the course of \cite{beckerkechris} that for $x_1,x_2\in X$, $x_1E^X_Gx_2$ if and only if 
there is some $g\in G$ with $g\cdot M_{x_1}=M_{x_2}$. 

At this point we may define $N_x$ be the expansion of $M_x$ obtained by incorperating all the relations of $M$. 
Since any $g\in G$ fixes $M_x$, we then obtain that that $N_{x_1}\cong N_{x_2}$ if and only if 
$\exists g\in G(g\cdot M_{x_1}=M_{x_2})$. We let ${\cal L}'$ be the language of these models $N_x$. 
Since $\{g\cdot N_x: g\in S_{\infty}\}$ is a Borel 
$S_{\infty}$ set, we may characterize it as the models of some $\sigma\in{\cal L}'_{\omega_1\omega}$.$\Box$\\

6.7. Theorem(Gao). Let $G=$Aut$(M)$ be a closed subgroup of $S_{\infty}$, with 
${\cal L}$ the language of $M$. 
Then $G$ is cli if and only if there is an ${\cal L}_{\omega_1\omega}$ elementary 
embedding $\pi: M\rightarrow M$ that is not onto.\\

6.8. Theorem(Knight). There is a countable model $M$ with language $\{<, f_0,f_1,...\}$, 
where 

(i) $<$ is a linear ordering on $M$; 

(ii) each $f_n$ is unary function; 

(iii) for each $y\in M$, $\{x\in M:x<y\}=\{f_n(y):n\in\omega\}$; and 

(iv) there is non-onto ${\cal L}_{\omega_1\omega}$ elementary 
embedding from $M$  to $M$.\\

6.9. Theorem. There is a Polish group $G$ such that:

(i) TVC$(G,\Ubf{\Sigma}^1_1)$ -- in the sense that if $A\subset X$ is $\Ubf{\Sigma}^1_1$, 
and $X$ is a Polish $G$-space, then either $|A/G|\leq\aleph_0$ or there is a perfect set 
$P\subset A$ of inequivalent reals; and 

(ii) $G$ is not cli. 

Proof. Let $G=$Aut$(M)$ for $M$ as in 6.8. 
$G$ is not cli by 6.8(iv) and 6.7. 
Suppose for a contradiction that TVC$(G,\Ubf{\Sigma}^1_1)$ fails, 
so let $X$ be a Polish $G$-space, $A\subset X$ 
with exactly $\aleph_1$ many orbits, and fix $\theta:X\rightarrow$ Mod$(\varphi)$ for some 
countable $\varphi\in{\cal L}_{\omega_1\omega}$ implying the Scott sentence of $M$, 
${\cal L}\supset\{<, f_0,f_1,...\}$. 

Thus as  in the proof of 5.3, 
${\Bbb P}_0=$Coll$(\omega,\kappa)$, $\kappa$ sufficiently big, 
there will some term $\sigma_0$ for ${\Bbb P}_0$ and $p_0\in{\Bbb P}_0$ 
with 
\[{\Bbb P}_0\Vdash \forall x\in {\Bbb V}(x\not\in[\sigma_0
[\dot{G}]]_G),\]
\[(p_0,p_0)\Vdash_{{\Bbb P}_0\times{\Bbb P_0}}\sigma_0[\dot{G_l}]E^X_G \sigma_0[\dot{G_r}].\] 
By applying this argument in again in ${\Bbb V}^{{\Bbb P_0}}$ we may find ${\cal P}_1$ and 
$\sigma_1$ such that 
\[{\Bbb P}_1\Vdash \forall x\in {\Bbb V}(x\not\in[\sigma_0
[\dot{G}]]_G),\]
\[(p_1,p_1)\Vdash_{{\Bbb P}_1\times{\Bbb P_1}}\sigma_1[\dot{G_l}]E^X_G \sigma_1[\dot{G_r}],\] 
\[{\Bbb P}_1\times{\Bbb P_0}\Vdash \sigma_0[\dot{G_l}]
\not\in[\sigma_0
[\dot{G_r}]]_G.\]

Continuing this transfinitely we may find $(\sigma_{\alpha}, {\Bbb P}_{\alpha}, p_{\alpha})$ for 
$\alpha$ an ordinal, such that for $\alpha\neq\beta$
\[(p_{\alpha}
,p_{\alpha})\Vdash_{{\Bbb P}_{\alpha}\times{\Bbb P}_{\alpha}}\sigma_{\alpha}[\dot{G_l}]E^X_G \sigma_{\alpha}[\dot{G_r}],\] 
\[{\Bbb P}_{\alpha}\times{\Bbb P}_{\beta}
\Vdash \sigma_{\alpha}[\dot{G_l}]
\not\in[\sigma_{\beta}
[\dot{G_r}]]_G.\]

Now for any such $\alpha$ and $[H_{\alpha}]\subset{\Bbb P}_{\alpha}$ ${\Bbb V}$-generic 
below $p_{\alpha}$, the equivalence class of $\sigma_{\alpha}[H_{\alpha}]$ does not 
depend on the choice of $H_{\alpha}$. Thus the isomorphism type of 
$\theta(\sigma_{\alpha}[H_{\alpha}])$ is independent of the choice of the 
generic, and thus so too the Scott sentence. Hence, as in $\S1$ of \cite{minimal},  
an induction on the set theoretical rank shows that the Scott sentence $\varphi_{\alpha}$ of 
$\theta(\sigma_{\alpha}[H_{\alpha}])$ exists in ${\Bbb V}$.

let $\gamma(\alpha)$ be the Scott height of any model of $\varphi_{\alpha}$ (where this model, 
as opposed to its Scott sentence,  
may only exist in a generic extension).  Let $A_{\alpha}$ be the 
collection of canonical $\gamma(\alpha)$ types realised by any such model. 
Note that the cardinality of $A_{\alpha}$ must be atleast that of 
$\gamma(\alpha)$, since for $\delta<\gamma(\alpha)$, and 
$N$ a model of $\varphi_{\alpha}$ appearing in some generic extension, there 
will be $\vec a, \vec b\in N^{<\omega}$ with $\varphi_{\delta}^{\vec a, N}=
\varphi_{\delta}^{\vec b, N}$ but 
$\varphi_{\delta+1}^{\vec a, N}\neq
\varphi_{\delta+1}^{\vec b, N}$. 
Note that $\gamma(\alpha)\rightarrow\infty$ as 
$\alpha\rightarrow\infty$, or else it would not possible for this ordinal sequence of 
Scott sentences to be non-repeating; and so $|A_{\alpha}|\rightarrow\infty$ as 
$\alpha\rightarrow\infty$. 

Now for any such $\alpha$ we can define a quasi-linear ordering on 
$A_{\alpha}$ by $\varphi'\leq\varphi''$ if and only if for any model $N$ of 
$\varphi(\alpha)$ and $\vec a, \vec b\in N$ with 
$\varphi_{\gamma(\alpha)}^{\vec a, N}=\varphi'$ and $\varphi_{\gamma(\alpha)}^{\vec b, N}=\varphi''$, 
for all $a_0\in\vec a$ there is some $b_0\in\vec b$ and some $c\in N$ with 
\[N\models c< b_0,\]
\[\varphi_{\gamma(\alpha)}^{a_0, N}=\varphi_{\gamma(\alpha)}^{c, N}.\]

Since each model of $\varphi$ is an expansion of the Knight model, it 
follows by 6.8(iii) that for each $\varphi'\in A(\alpha)$ there are at most 
countably many $\varphi''\leq\varphi'$. Now taking some 
$A(\alpha)$ with size bigger than $\aleph_1$ we have a contradiction.$\Box$\\

A positive answer to the next question would help clarify 6.9.\\

6.10 Question. If $G=$Aut$(M)$, for $M$ a countable model, and TVC$(G,\Ubf{\Sigma}^1_1)$ fails, then must 
$M$ have a model of size $2^{\aleph_0}$?



6363 MSB

Mathematics

UCLA

CA90095-1555

greg@math.ucla.edu

\end{document}